\documentclass[12pt]{article}

\usepackage{latexsym,amsmath}
\usepackage{amssymb}
\usepackage{url}
\usepackage{graphicx}
\usepackage{microtype}

\usepackage[margin=1.18in]{geometry}

\newtheorem{thm}{Theorem}[section]

\newtheorem{lemma}[thm]{Lemma}

\newtheorem{conj}[thm]{Conjecture}

\newcommand{\qed}[0]{{\hspace*{\fill}\mbox{$\Box$}}}

\begin{document}

\renewcommand{\thefootnote}{\fnsymbol{footnote}}

\title{Maximizing $H$-colorings of a regular graph}

\author{David Galvin\thanks{Department of Mathematics,
University of Notre Dame, 255 Hurley Hall, Notre Dame IN
46556; dgalvin1@nd.edu. Research supported in part by National Security Agency grant H98230-10-1-0364.}}

\maketitle

\begin{abstract}
For graphs $G$ and $H$, a {\em homomorphism} from $G$ to $H$, or {\em $H$-coloring} of $G$, is an adjacency preserving map from the vertex set of $G$ to the vertex set of $H$. Our concern in this paper is the maximum number of $H$-colorings admitted by an $n$-vertex, $d$-regular graph, for each $H$. Specifically, writing ${\rm hom}(G,H)$ for the number of $H$-colorings admitted by $G$, we conjecture that for any simple finite graph $H$ (perhaps with loops) and any simple finite $n$-vertex, $d$-regular, loopless graph $G$ we have
$$
{\rm hom}(G,H) \leq \max\left\{{\rm hom}(K_{d,d},H)^{\frac{n}{2d}}, {\rm hom}(K_{d+1},H)^{\frac{n}{d+1}}\right\}
$$
where $K_{d,d}$ is the complete bipartite graph with $d$ vertices in each partition class, and $K_{d+1}$ is the complete graph on $d+1$ vertices.

Results of Zhao confirm this conjecture for some choices of $H$ for which the maximum is achieved by ${\rm hom}(K_{d,d},H)^{n/2d}$. Here we exhibit for the first time infinitely many non-trivial triples $(n,d,H)$ for which the conjecture is true and for which the maximum is achieved by ${\rm hom}(K_{d+1},H)^{n/(d+1)}$.

We also give sharp estimates for ${\rm hom}(K_{d,d},H)$ and ${\rm hom}(K_{d+1},H)$ in terms of some structural parameters of $H$. This allows us to characterize those $H$ for which ${\rm hom}(K_{d,d},H)^{1/2d}$ is eventually (for all sufficiently large $d$) larger than ${\rm hom}(K_{d+1},H)^{1/(d+1)}$ and those for which it is eventually smaller, and to show that this dichotomy covers all non-trivial $H$. Our estimates also allow us to obtain asymptotic evidence for the conjecture in the following form. For fixed $H$, for all $d$-regular $G$ we have
$$
{\rm hom}(G,H)^{\frac{1}{|V(G)|}} \leq (1+o(1))\max\left\{{\rm hom}(K_{d,d},H)^{\frac{1}{2d}}, {\rm hom}(K_{d+1},H)^{\frac{1}{d+1}}\right\}
$$
where $o(1)\rightarrow 0$ as $d \rightarrow \infty$. More precise results are obtained in some special cases.
\end{abstract}

\section{Introduction}

A {\em homomorphism} from a simple loopless finite graph $G=(V(G),E(G))$ to a simple finite graph $H=(V(H),E(H))$ (perhaps with loops) is an adjacency preserving function from the vertices of $G$ to the vertices of $H$. We write
$$
{\rm Hom}(G,H) = \{f:V(G)\rightarrow V(H) | uv \in E(G) \Rightarrow f(u)f(v) \in E(H)\}
$$
for the set of all homomorphisms from $G$ to $H$, and ${\rm hom}(G,H)$ for $|{\rm Hom}(G,H)|$. (All graphs in this paper will be simple and finite. Those denoted by $G$ will always be loopless, but loops will be allowed in those denoted by $H$.)

Graph homomorphisms generalize some important notions in graph theory. When $H=H_{\rm ind}$ consists of a single edge with a loop at one vertex, ${\rm Hom}(G,H)$ can be identified with the set of independent sets of $G$ (sets of vertices no two of which are adjacent), via the inverse image of the unlooped vertex. When $H=K_q$ (the complete graph on $q$ vertices), ${\rm Hom}(G,H)$ can be identified with the set of proper $q$-colorings of $G$ (assignments of labels to the vertices of $G$, from a palette of $q$ possible labels, with the property that adjacent vertices receive different labels). Motivated by the latter example, an element of ${\rm Hom}(G,H)$ is often referred to as an {\em $H$-coloring} of $G$, and this is the terminology that we will use throughout this paper.

In statistical physics, $H$-colorings arise as configurations in {\em hard-constraint spin systems}. Here we think of the vertices of $G$ as locations, each one occupied by a particle having one of a set of spins or colors indexed by vertices of the {\em constraint graph} $H$. The edges of $G$ encode pairs of spins that are bonded (for example by spatial proximity) and the occupation rule is that the spins appearing at a pair of bonded locations must
be adjacent in $H$. A valid configuration of spins on $G$ is thus exactly an $H$-coloring of $G$. In the language of statistical physics, independent sets are configurations in the {\em hard-core gas model} (a model of the occupation of space by a gas with massive particles; the unlooped vertex in $H_{\rm ind}$ represents the presence of a particle and the looped vertex represents the absence of a particle), while proper $q$-colorings are configurations in the {\em zero-temperature $q$-state antiferromagnetic Potts model}. Another constraint graph that appears frequently in statistical physics and will play a role in the present paper is the {\em Widom-Rowlinson} graph $H_{\rm WR}$. This is the completely looped path on three vertices. If the end vertices of the path are considered to be two different types of particles, with the middle vertex representing empty space, then the Widom-Rowlinson constraint graph models the occupation of space by two mutually repulsive particles.

\medskip

A natural question to ask concerning $H$-colorings is, which graphs $G$ in a given family maximize (or minimize) ${\rm hom}(G,H)$? For the family of graphs with a fixed number of vertices and edges, this question for $H=K_q$ was first asked independently by Wilf and Linial around 1986, and although it has not yet been answered completely, significant progress has been made (see \cite{LohPikhurkoSudakov} and the references therein). For the same family (fixed number of vertices and edges) the question was raised and essentially completely answered for $H=H_{\rm ind}$ and $H=H_{\rm WR}$ by Cutler and Radcliffe in \cite{CutlerRadcliffe2}.

For $H=H_{\rm ind}$ and the family of regular graphs, the question was first raised by Granville in the context of combinatorial group theory (see \cite{Alon}) and an approximate answer was given by Alon \cite{Alon}. A more complete answer was given by Kahn \cite{Kahn}, who used entropy methods to show that if $G$ is an $n$-vertex, $d$-regular {\em bipartite} graph then ${\rm hom}(G,H_{\rm ind}) \leq {\rm hom}(K_{d,d},H_{\rm ind})^{n/2d}$, where $K_{d,d}$ is the complete bipartite graph with $d$ vertices in each class. Note that this bound is tight, being achieved in the case $2d|n$ by the graph consisting of the disjoint union of $n/2d$ copies of $K_{d,d}$. Also using entropy, Galvin and Tetali \cite{GalvinTetali-weighted} generalized Kahn's result to arbitrary $H$, showing that for all $n$-vertex, $d$-regular, bipartite $G$ we have
\begin{equation} \label{GTbound}
{\rm hom}(G,H) \leq {\rm hom}(K_{d,d},H)^\frac{n}{2d}.
\end{equation}

Both Alon and Kahn conjectured that in the case $H=H_{\rm ind}$, (\ref{GTbound}) remains true if $G$ is allowed to vary over all $d$-regular graphs (not necessarily bipartite); in other words, the number of independent sets admitted by an $n$-vertex, $d$-regular graph is never more than the number admitted by $K_{d,d}$ raised to the power $n/2d$. This conjecture was resolved by Zhao \cite{Zhao}, who used an elegant argument to deduce the general result from the bipartite case.

Galvin and Tetali proposed that (\ref{GTbound}) should hold for all $H$ even if $G$ is not assumed to be bipartite. As we will presently discuss, Zhao \cite{Zhao, Zhao2} and Lazebnik (personal communication) showed that this stronger conjecture is true for $H$ in various classes, but it is false in general: Cutler, Radcliffe (personal communications) and Zhao \cite{Zhao2} have all provided counter-examples, the simplest being the disjoint union of two loops.

A common feature of all the counter-examples is that they are graphs $H$ for which
$$
{\rm hom}(K_{d,d},H)^{\frac{1}{2d}} < {\rm hom}(K_{d+1},H)^{\frac{1}{d+1}}
$$
for some suitably large $d$, where $K_{d+1}$ is the complete graph on $d+1$ vertices. So for these $H$ and for an $n$ which is a multiple of both $2d$ and $d+1$, there are more $H$-colorings of the graph consisting of $n/(d+1)$ copies of $K_{d+1}$ than there are of the graph consisting of $n/2d$ copies of $K_{d,d}$. No example is known of an $H$ and an $n$-vertex, $d$-regular $G$ for which
$$
{\rm hom}(G,H)>\max\left\{{\rm hom}(K_{d,d},H)^{\frac{n}{2d}}, {\rm hom}(K_{d+1},H)^{\frac{n}{d+1}}\right\},
$$
and is tempting to propose the following conjecture.
\begin{conj} \label{conj-GTcorr}
Fix $n$, $d$ and $H$. For all $n$-vertex $d$-regular graphs $G$ we have
$$
{\rm hom}(G,H) \leq \max\left\{{\rm hom}(K_{d,d},H)^{\frac{n}{2d}}, {\rm hom}(K_{d+1},H)^{\frac{n}{d+1}}\right\}.
$$
\end{conj}
Note that because ${\rm hom}(G_1 \cup G_2,H)={\rm hom}(G_1,H){\rm hom}(G_2,H)$ (where $G_1 \cup G_2$ is the union of disjoint copies of $G_1$ and $G_2$), Conjecture \ref{conj-GTcorr} for arbitrary $G$ follows from the special case of connected $G$.

To avoid trivialities we assume throughout that $d \geq 1$. Given this assumption, we may (and will) also assume that $H$ has no isolated vertices, since if $G$ is any graph with minimum degree at least one, $H$ is any graph, and $H'$ is obtained from $H$ by adding some isolated vertices, then we have ${\rm hom}(G,H')={\rm hom}(G,H)$ (and if $H$ is itself an isolated vertex, then ${\rm hom}(G,H)=0$ for all $G$).

If $H$ is a complete looped graph on $k$ vertices then for all $n$-vertex $G$ we have
$$
{\rm hom}(G,H) = {\rm hom}(K_{d,d},H)^{\frac{n}{2d}} =  {\rm hom}(K_{d+1},H)^{\frac{n}{d+1}} = k^n
$$
and so for these $H$ Conjecture \ref{conj-GTcorr} is trivially true for all $n$ and $d$. We will see in Section \ref{sec-recipe} that this is the only circumstance under which there is equality between ${\rm hom}(K_{d,d},H)^{1/2d}$ and ${\rm hom}(K_{d+1},H)^{1/(d+1)}$ for all $d$; in fact, for each $H$ that is not a completely looped graph we have either ${\rm hom}(K_{d+1},H)^{1/(d+1)} < {\rm hom}(K_{d,d},H)^{1/2d}$ for all large $d$ (we say that these $H$ are of {\em complete bipartite type}) or we have ${\rm hom}(K_{d,d},H)^{1/2d} < {\rm hom}(K_{d+1},H)^{1/(d+1)}$ for all large $d$ ({\em complete type}). Moreover, there is a simple procedure that allows us determine the type of a given $H$; see Theorem \ref{thm-char}.

\medskip

There are two other trivial instances of Conjecture \ref{conj-GTcorr}. First, if $H$ is bipartite then for non-bipartite $G$ we have ${\rm hom}(G,H)=0$ and so in this case the conjecture follows for all $n$ and $d$ from (\ref{GTbound}) (and $H$ is of complete bipartite type). Second, if $H$ is the disjoint union of $k>2$ loops then ${\rm hom}(G,H)=k^{{\rm comp}(G)}$ (where ${\rm comp}(G)$ is the number of components of $G$). This is in turn at most $k^{n/(d+1)}={\rm hom}(K_{d+1},H)^{n/(d+1)}$ (since each component of a $d$-regular $G$ must have at least $d+1$ vertices). So for these $H$ (which are of complete type) the conjecture is also true for all $n$ and $d$.

There are also some non-trivial examples of $H$ for which the conjecture is known to be true. Given a graph $H$, construct $H^{\rm bst}$ on vertex set $V(H) \times V(H)$ as follows: join $(u,v)$ and $(u',v')$ if and only if i) $u$ and $u'$ are joined in $H$, ii) $v$ and $v'$ are joined in $H$, and iii) either $u$ and $v'$ or $u'$ and $v$ are not joined in $H$. Zhao \cite{Zhao2} showed that if $H^{\rm bst}$ is bipartite then for each non-bipartite $d$-regular $G$ we have
${\rm hom}(G,H) < {\rm hom}(K_{d,d},H)^{n/2d}$,
and so for these $H$ (which are of complete bipartite type) we have Conjecture \ref{conj-GTcorr} for all $n$ and $d$ by (\ref{GTbound}). An example of an $H$ for which $H^{\rm bst}$ is bipartite is $H_{\rm ind}$; Zhao had earlier \cite{Zhao} dealt with this special case.

Zhao observed that all threshold-type graphs $H$ have the property that $H^{\rm bst}$ is bipartite, where a {\em threshold-type graph} is one for which there exists some assignment $\alpha: V(H) \rightarrow {\mathbb R}$ and some $t \in {\mathbb R}$ such that $uv \in E(H)$ if and only if $\alpha(u) + \alpha(v) \leq t$ for every (not necessarily distinct) $u, v \in V(H)$. Examples of threshold-type graphs include $H_{\rm ind}$ and more generally $H(k)$, the $k$-state hard-core constraint graph (which will be discussed after the proof of Theorem \ref{thm-char}).

Zhao \cite{Zhao2} also looked at the case $H=K_q$ (proper $q$-coloring), and showed that if $G$ is an $n$-vertex, $d$-regular graph and $q$ is sufficiently large (specifically, $q>(2n)^{2n+2}$) then ${\rm hom}(G,K_q) \leq {\rm hom}(K_{d,d},K_q)^{n/2d}$. Applying this result for $G=K_{d+1}$ we see that this verifies Conjecture \ref{conj-GTcorr} for all triples $(n,d,K_q)$ with $q>(2n)^{2n+2}$.

A method used by Lazebnik on a related problem can be used to significantly improve the bound on $q$ in the case when $2d|n$.
\begin{thm} \label{thm_Laz}
Fix $n$ and $d$ with $2d|n$. Let $G$ be an $n$-vertex, $d$-regular graph. For $q > 2{nd/2 \choose 4}$ we have
$$
{\rm hom}(G,K_q) \leq {\rm hom}(K_{d,d},K_q)^\frac{n}{2d},
$$
with equality if and only if $G$ is the disjoint union of $n/2d$ copies of $K_{d,d}$.
\end{thm}
Lazebnik's approach uses Whitney's broken circuit theorem \cite{Whitney}, which provides a combinatorial interpretation of the coefficients of the chromatic polynomial (we give the details in Section \ref{sec-genincexc}). A related but more elementary approach allows us to tackle a more general class of graphs $H$. The idea is based on the following lemma, whose proof is a direct application of the principle of inclusion-exclusion.
\begin{lemma} \label{inq-incexc}
For any simple, finite $H$ and any simple, finite, loopless $G$, we have
$$
{\rm hom}(G,H) = \sum_{S \subseteq E(G)} (-1)^{|S|}{\rm hom}(G(S),H^c)|V(H)|^{n-v(S)}
$$
where $v(S)$ is the number of vertices spanned by the edge set $S$, $G(S)$ is the subgraph of $G$ spanned by $S$, and $H^c$ is the complement of $H$ (with all isolated vertices removed).
\end{lemma}
If $H^c$ is fairly simple, then it may be possible to understand the terms ${\rm hom}(G(S),H^c)$ and so to understand the behavior of ${\rm hom}(G,H)$ for large $|V(H)|$ well enough to obtain a precise result. In Section \ref{sec-genincexc} we give an example of this approach, to prove the following theorem.
\begin{thm} \label{thm_Hqell}
Fix $n$ and $d$ with $2d|n$. Let $G$ be an $n$-vertex, $d$-regular graph. Let $H_q^\ell$ be the complete looped graph on $q$ vertices with $\ell \geq 1$ loops deleted. Then for $q > \exp_2\left\{nd/2\right\}$ we have
$$
{\rm hom}(G,H_q^\ell) \leq {\rm hom}(K_{d,d},H_q^\ell)^\frac{n}{2d},
$$
with equality if and only if $G$ is the disjoint union of $n/2d$ copies of $K_{d,d}$.
\end{thm}
This verifies Conjecture \ref{conj-GTcorr} for triples $(n,d,H_q^\ell)$ with $q > \exp_2\left\{nd/2\right\}$ and $2d(d+1)|n$.

As a special case of Theorem \ref{thm_Hqell} (the case $\ell=q$), we recover Theorem \ref{thm_Laz} (on proper $q$-colorings), with a bound on $q$ that is worse than that of Theorem \ref{thm_Laz} but better than Zhao's bound from \cite{Zhao2}. Another interesting special case is $\ell=1$. The number of elements of ${\rm hom}(G,H_q^1)$ in which $k$ vertices of $G$ get mapped to the unlooped vertex of $H_q^1$ is easily seen to be $(q-1)^{n-k}$ times the number of independent sets in $G$ of size $k$. So for fixed $n$ and $d$ with $2d|n$ and for $q>\exp_2\left\{nd/2\right\}$ we have that for all $n$-vertex, $d$-regular $G$,
\begin{eqnarray*}
\sum_{k=0}^n i_k(G)(q-1)^{n-k} & = & {\rm hom}(G,H_q^1) \\
& \leq & {\rm hom}\left(\frac{n}{2d}K_{d,d},H_q^1\right) \\
& = & \sum_{k=0}^n i_k\left(\frac{n}{2d}K_{d,d}\right)(q-1)^{n-k}
\end{eqnarray*}
where $\frac{n}{2d}K_{d,d}$ is the disjoint union of $n/2d$ copies of $K_{d,d}$ and $i_k(G)$ is the number of independent sets in $G$ of size $k$. Scaling by $(q-1)^n$ we get
$$
\sum_{k=0}^n i_k(G)\left(\frac{1}{q-1}\right)^k \leq \sum_{k=0}^n i_k\left(\frac{n}{2d}K_{d,d}\right)\left(\frac{1}{q-1}\right)^k,
$$
with equality if and only if $G=\frac{n}{2d}K_{d,d}$. This is a special case of a result Zhao \cite{Zhao}, building on results of Kahn \cite{Kahn,Kahn3} and Galvin and Tetali \cite{GalvinTetali-weighted}, that states that for all $\lambda > 0$ we have
$$
\sum_{k=0}^n i_k(G)\lambda^k \leq \sum_{k=0}^n i_k\left(\frac{n}{2d}K_{d,d}\right)\lambda^k.
$$
The proof of Zhao's result relies on the entropy method. To the best of our knowledge this is the first proof of any special case that avoids the use of entropy.

\medskip

There are no non-trivial examples of triples $(n,d,H)$ for which Conjecture \ref{conj-GTcorr} has been proven and for which the maximum is achieved by ${\rm hom}(K_{d+1},H)^{n/(d+1)}$. Using the inclusion-exclusion approach outlined above, we provide the first such examples here. Let $H_q$ be the complete looped graph on $q$ vertices with one edge (not a loop) removed. We may think of this as a ``$q$-state Widom-Rowlinson'' constraint graph (in the case $q=3$, it is exactly $H_{\rm WR}$ introduced earlier). In Section \ref{sec-WRgen} we prove that if $n$, $d$ and $q$ satisfy $(d+1)|n$ and $q > \exp_2\left\{nd/2+n/2-1\right\}$ then for any $n$-vertex, $d$-regular graph $G$ we have
$$
{\rm hom}(G,H_q) \leq {\rm hom}(K_{d+1},H_q)^\frac{n}{d+1},
$$
with equality if and only if $G$ is the disjoint union of $n/(d+1)$ copies of $K_{d+1}$. More generally, we have the following.
\begin{thm} \label{thm-qstateWR}
Let $\ell \geq 1$, $n$ and $d$ be given, with $(d+1)|n$. There is $q_0=q_0(n,d,\ell)$ such that if $q\geq q_0$ and $H$ is obtained from the complete looped graph on $q$ vertices by the deletion of a collection of no more than $\ell$ edges spanning a disjoint union of complete bipartite subgraphs, then for any $n$-vertex, $d$-regular graph $G$ we have
$$
{\rm hom}(G,H) \leq {\rm hom}(K_{d+1},H)^\frac{n}{d+1},
$$
with equality if and only if $G$ is the disjoint union of $n/(d+1)$ copies of $K_{d+1}$.
\end{thm}
This verifies Conjecture \ref{conj-GTcorr} for all triples $(n,d,H)$ with $2d(d+1)|n$ and $H$ as described above.

\medskip

In Section \ref{sec-recipe} we obtain sharp estimates on both ${\rm hom}(K_{d,d},H)$ and ${\rm hom}(K_{d+1},H)$ in terms of some structural parameters of $H$. We use these estimates in Section \ref{sec-asymptotics} to examine the asymptotics (in $d$) of the quantities in Conjecture \ref{conj-GTcorr}, and give the following asymptotic verification.
\begin{thm} \label{thm-asymptotic}
Fix $H$. For all $d$-regular $G$,
$$
{\rm hom}(G,H)^{\frac{1}{|V(G)|}} \leq (1+o(1))\max\left\{{\rm hom}(K_{d,d},H)^{\frac{1}{2d}}, {\rm hom}(K_{d+1},H)^{\frac{1}{d+1}}\right\}
$$
where $o(1)\rightarrow 0$ as $d \rightarrow \infty$.
\end{thm}
For some $H$ we can obtain more precise statements; the details are in Section \ref{sec-asymptotics}.

The estimates of Section \ref{sec-asymptotics} also allow us to obtain some exact results. For example we can prove that if ${\rm hom}(G,H_{\rm WR})$ exceeds its conjectured upper bound, then $G$ must be close to being a disjoint union of $n/(d+1)$ copies of $K_{d+1}$.
\begin{thm} \label{thm-WR}
There is a constant $C>0$ such that for all sufficiently large $d$, we have the following. If $G$ is an $n$-vertex, $d$-regular graph with an independent set of size at least $Cn/d$, then
$$
{\rm hom}(G,H_{\rm WR}) < {\rm hom}(K_{d+1},H_{\rm WR})^\frac{n}{d+1}.
$$
\end{thm}

To justify the assertion that those $G$ not covered by Theorem \ref{thm-WR} are close to being a disjoint union of $n/(d+1)$ copies of $K_{d+1}$, we use two metrics -- the size of the maximal independent set in $G$, and the average local density.

All $n$ vertex, $d$-regular graphs with $(d+1)|n$ have maximum independent set size at least $n/(d+1)$, with $n/(d+1)$ being achieved uniquely by the disjoint union of $n/(d+1)$ copies of $K_{d+1}$ (this is Tur\'an's theorem). Moreover, for $3 \leq d\leq n^{1-\Omega(1)}$ the largest independent set in a uniformly chosen $d$-regular graph on $n$ vertices has size concentrated close to $2n\log d/d$ (see \cite{CooperFriezeReedRiordan}). So for $d$ in this range Theorem \ref{thm-WR} says that Conjecture \ref{conj-GTcorr} is almost surely true for $H=H_{\rm WR}$; and to prove the conjecture for all $G$ it (just) requires bringing the constant $C$ down to around $1$.

A result of Shearer \cite{Shearer} shows that if an $n$-vertex, $d$-regular graph $G$ has at most $m$ triangles then it has an independent set of size at least $(4n/78d)\log (d^2n/m)$. It follows that there is a constant $c$ (depending on the $C$ from Theorem \ref{thm-WR}) such that for all sufficiently large $d$, if $n$-vertex, $d$-regular $G$ has fewer than $d^2n/c$ triangles then Conjecture \ref{conj-GTcorr} (in the case $H=H_{\rm WR}$) is true for $G$. This means that the graphs for which the conjecture remains open must be on average locally dense: a pair of adjacent vertices must have on average $\Omega(d)$ common neighbors. Note that a pair of adjacent vertices in the disjoint union of $n/(d+1)$ copies of $K_{d+1}$ always has $d-1$ common neighbors, and this is the maximum possible average local density.

We can extend Theorem \ref{thm-WR} to an infinite family of graphs; we postpone the details to Section \ref{sec-asymptotics}, as the statement requires notation from Section \ref{sec-recipe}.

\medskip

The remainder of this paper is laid out as follows. In Section \ref{sec-genincexc} we use the inclusion-exclusion approach to prove Theorems \ref{thm_Laz} and \ref{thm_Hqell}. In Section \ref{sec-WRgen} we use the same approach to prove Theorem \ref{thm-qstateWR}. In Section \ref{sec-recipe} we obtain sharp estimates on ${\rm hom}(K_{d,d},H)$ and ${\rm hom}(K_{d+1},H)$, and use these to characterize those $H$ for which ${\rm hom}(K_{d,d},H)^{1/2d}$ is eventually larger (or smaller) than ${\rm hom}(K_{d+1},H)^{1/(d+1)}$ (Theorem \ref{thm-char}). In Section \ref{sec-asymptotics} we use our estimates from Section \ref{sec-recipe} to prove Theorem \ref{thm-asymptotic}, along with other more precise estimates, including Theorem \ref{thm-WR-gen} (a generalization of Theorem \ref{thm-WR}).

\section{Proofs of Theorems \ref{thm_Laz} and \ref{thm_Hqell}} \label{sec-genincexc}

In what follows we denote by $\frac{n}{2d}K_{d,d}$ the graph consisting of the disjoint union of $n/2d$ copies of $K_{d,d}$.
We need the following lemma. Here we denote by $c_4(G)$ the number of (not necessarily induced) cycles on four vertices contained in $G$.
\begin{lemma} \label{lemma-c4}
Let $G$ be a $d$-regular triangle-free graph on $n$ vertices with $2d|n$. We have
$$
c_4(G) \leq c_4\left(\frac{n}{2d}K_{d,d}\right)
$$
with equality if and only if $G=\frac{n}{2d}K_{d,d}$.
\end{lemma}

\medskip

\noindent {\em Proof}: The number of $4$-cycles passing through a particular edge $xy$ is at most $(d-1)^2$. This bound is achieved if all $d-1$ neighbors of $x$ (other than $y$) are joined to all $d-1$ neighbors of $y$ (other than $x$), that is, if the component of the edge $xy$ is $K_{d,d}$. Since $c_4(G)$ is one quarter the sum over all edges of the number of $4$-cycles passing through that edge, it follows that the total number of $4$-cycles in $G$ is at most $c_4(\frac{n}{2d}K_{d,d})$, and that $\frac{n}{2d}K_{d,d}$ is the unique graph that achieves the maximum.
\qed

\medskip

\noindent {\em Proof of Theorem \ref{thm_Hqell}}: Specifying Lemma \ref{inq-incexc} to $H_q^\ell$ we get
\begin{equation} \label{inq-incexc2}
{\rm hom}(G,H_q^\ell) = \sum_{S \subseteq E(G)} (-1)^{|S|}\ell^{c(S)}q^{n-v(S)}
\end{equation}
where $c(S)$ and $v(S)$ are the number of components and vertices of the subgraph spanned by $S$.

We examine (\ref{inq-incexc2}) first in the case where $G$ has at least one triangle.
The contribution to (\ref{inq-incexc2}) from $|S|=0$, $1$ and $2$ is independent of $G$ (depending only on $n$ and $d$). For example, the contribution from those $S$ with $|S|=2$ that span a path on three vertices is $n{d \choose 2}\ell q^{n-3}$ (for each vertex $v$ there are ${d \choose 2}$ choices of such a path with $v$ as the center vertex), while the contribution from those $S$ with $|S|=2$ that do not span a path is $({nd/2 \choose 2}-n{d \choose 2})\ell^2 q^{n-4}$
(note that $nd/2$ is the number of edges in $G$). We lower bound ${\rm hom}(\frac{n}{2d}K_{d,d},H_q^\ell)$ by taking the terms on the right-hand side of (\ref{inq-incexc2}) corresponding to $|S|=0$, $1$ and $2$, and the terms corresponding to $|S| \geq 3$ and odd. We upper bound ${\rm hom}(G,H_q^\ell)$ by taking the terms on the right-hand side of (\ref{inq-incexc2}) corresponding to $|S|=0$, $1$ and $2$, the terms corresponding to $|S| \geq 3$ and even, and the contribution from a single $S$ that spans a triangle. This yields
\begin{eqnarray*}
{\rm hom}\left(\frac{n}{2d}K_{d,d},H_q^\ell\right) - {\rm hom}(G,H_q^\ell) & \geq & \ell q^{n-3} - \sum_{S \subseteq \tilde{E}} \ell^{c(S)}q^{n-v(S)} \\
& \geq & \ell q^{n-3} - 2^{nd/2} \max\left\{\ell^{c(S)}q^{n-v(S)}\right\}
\end{eqnarray*}
where $\tilde{E}$ is taken to be the edge set of $\frac{n}{2d}K_{d,d}$ if $|S|$ is odd, and of $G$ if $|S|$ is even, and the sum and maximum are both over those $S$ with $|S| \geq 3$. To understand the maximum, first note that $c(S)\leq \lfloor v(S)/2 \rfloor$ (since each component must have at least $2$ vertices). For $v(S)\geq 5$ we have $v(S)-\lfloor v(S)/2 \rfloor \geq 3$ and so
$$
\ell^{c(S)}q^{n-v(S)} \leq \ell^{\lfloor v(S)/2 \rfloor}q^{n-v(S)} \leq \ell q^{n-v(S) + \lfloor v(S)/2 \rfloor -1} \leq \ell q^{n-4}.
$$
Also, every $S$ with $v(S)=4$ and $|S|\geq 3$ has $c(S)=1$ and so $\ell^{c(S)}q^{n-v(S)} = \ell q^{n-4}$. This covers all possible $S$ (note that we do not consider $|S|=3$ spanning a triangle here, since $\frac{n}{2d}K_{d,d}$ is triangle free). It follows that
$$
{\rm hom}\left(\frac{n}{2d}K_{d,d},H_q^\ell\right) - {\rm hom}(G,H_q^\ell) \geq  \ell q^{n-3} - 2^{nd/2}\ell q^{n-4},
$$
which is strictly positive for all $q > 2^{nd/2}$.

Thus we may assume that $G$ is triangle-free, and different from $\frac{n}{2d}K_{d,d}$. Within the class of triangle-free graphs, the contribution to (\ref{inq-incexc2}) from $|S|=3$ is independent of $G$
(the number of subsets of three edges that span a star on four vertices, a path on four vertices, and the two-component graph consisting of a path on three vertices together with a edge, are all easily seen to be independent of $G$; the count for the one remaining subgraph, three disconnected edges, must therefore also be independent of $G$).
As before, the contributions from $|S|=0$, $1$ and $2$ are also independent of $G$. By Lemma \ref{lemma-c4} $G$ has fewer $4$-cycles than $\frac{n}{2d}K_{d,d}$ and so by the same process as in the case where $G$ has a triangle, we have
$$
{\rm hom}\left(\frac{n}{2d}K_{d,d},H_q^\ell\right) - {\rm hom}(G,H_q^\ell) \geq \ell q^{n-4} - 2^{nd/2} \max\left\{\ell^{c(S)}q^{n-v(S)}\right\}
$$
where now the maximum is over those $S$ with $|S| \geq 4$ and $S$ not spanning a $4$-cycle. For $v(S)\geq 7$ we have $v(S)-\lfloor v(S)/2 \rfloor \geq 4$, so reasoning as before we get that for these $S$ the quantity being maximized is at most $\ell q^{n-5}$. Those $S$ with $v(S)=6$ and $c(S)=3$ must span three edges, and so have $|S|=3$ and are not being considered. For all remaining $S$ with $v(S)=6$ we have $v(S)-c(S)\geq 4$ and so again the quantity being maximized is at most $\ell q^{n-5}$. For those $S$ with $v(S)=5$ and $c(S)=1$ we have $v(S)-c(S)= 4$ and so again the quantity being maximized is at most $\ell q^{n-5}$. In a triangle-free $G$ there are no $S$'s with $|S|\geq 4$, $v(S)=5$ and $c(S)>1$, and no $S$'s with $v(S)=4$ and $|S|\geq 4$ except those $S$ spanning a $4$-cycle. It follows that we have
$$
{\rm hom}\left(\frac{n}{2d}K_{d,d},H_q^\ell\right) - {\rm hom}(G,H_q^\ell) \geq  \ell q^{n-4} - 2^{nd/2}\ell q^{n-5},
$$
which is strictly positive for all $q > 2^{nd/2}$; the theorem is now proved. \qed

\medskip

\noindent {\em Proof of Theorem \ref{thm_Laz}}: For $d=1$ the result is trivial, as it is for $d=2$ and $n=4$. So we may assume that $nd \geq 16$.

We begin by recalling that for all graphs $G$, the quantity ${\rm hom}(G,K_q)$, viewed as a function of $q$, turns out to be a polynomial in $q$, the {\em chromatic polynomial}, which may be expressed as
\begin{equation} \label{Whitney}
{\rm hom}(G,K_q) = q^n + \sum_{i = 1}^{n-1} (-1)^i a_i(G)q^{n-i}
\end{equation}
where $n$ is the number of vertices in $G$, and the $a_i$'s are non-negative integers.  Whitney \cite{Whitney} gave a combinatorial interpretation of the $a_i$'s, in terms of broken circuits. A {\em broken circuit} in a graph is obtained from the edge set of a cycle by removing the maximum edge  of the cycle (with respect to some fixed
linear ordering of the edges); Whitney showed that $a_i$ is the number of subsets of size $i$ of the edge set of $G$ that do not contain a broken circuit as a subset.

Using (\ref{Whitney}) we can easily compute the first few coefficients of ${\rm hom}(G,K_q)$. Since a broken circuit must have size at least $2$, we have $a_1(G) = nd/2$ (the size of the edge set of $G$). Each triangle in $G$ gives rise to a pair of edges that may not be counted in the calculation of $a_2(G)$ and so
$a_2(G) = {nd/2 \choose 2} - c_3(G)$
where $c_3(G)$ is the number of triangles of $G$. In general $a_3(G)$ is a little harder to compute, except in the special case where $c_3(G)=0$. In this case each $4$-cycle in $G$ gives rise to a triple of edges that may not be counted in the calculation of $a_2(G)$ and so
$a_3(G) = {nd/2 \choose 3} - c_4(G)$
where $c_4(G)$ is the number of $4$-cycles of $G$. Directly from the definition we also have the general bound
$0 \leq a_i(G) \leq {nd/2 \choose i}$
valid for all $1 \leq i \leq n-1$.

For $G$ with $c_3(G) > 0$, all this together gives
$$
{\rm hom}\left(\frac{n}{2d}K_{d,d},K_q\right) - {\rm hom}(G,K_q) \geq q^{n-2} - \sum_{i=3}^{n-1} {nd/2 \choose i}q^{n-i}.
$$
Now for all $i=3, \ldots, n-2$ we have
$$
{nd/2 \choose i}q^{n-i} \geq 2{nd/2 \choose i+1}q^{n-(i+1)}
$$
as long as $q \geq (nd/2-3)/2$, so in this range we have
$$
\sum_{i=3}^{n-1} {nd/2 \choose i}q^{n-i} \leq 2{nd/2 \choose 3}q^{n-3}.
$$
It follows that ${\rm hom}(\frac{n}{2d}K_{d,d},K_q) - {\rm hom}(G,K_q) > 0$ as long as $q>2{nd/2 \choose 3}$.

For those $G$ with $c_3(G)=0$ we may assume that $c_4(G) < c_4(\frac{n}{2d}K_{d,d})$, since otherwise by Lemma \ref{lemma-c4} we have $G=\frac{n}{2d}K_{d,d}$ and the result is trivial. Now we have
$$
{\rm hom}\left(\frac{n}{2d}K_{d,d},K_q\right) - {\rm hom}(G,K_q) \geq q^{n-3} - \sum_{i=4}^{n-1} {nd/2 \choose i}q^{n-i},
$$
which by the same reasoning as before is strictly positive as long as $q>2{nd/2 \choose 4} \geq 2{nd/2 \choose 3}$ (the second inequality valid since $nd \geq 16$).

We conclude that for $q>2{nd/2 \choose 4}$ we have ${\rm hom}(G,K_q) \leq {\rm hom}(\frac{n}{2d}K_{d,d},K_q)$ for all $n$-vertex, $d$-regular $G$. \qed

\section{Proof of Theorem \ref{thm-qstateWR}} \label{sec-WRgen}

In what follows we denote by $\frac{n}{d+1}K_{d+1}$ the graph consisting of the disjoint union of $n/(d+1)$ copies of $K_{d+1}$.
We need the following lemma. Here we denote by $p_4(G)$ the number of (not necessarily induced) paths on four vertices in $G$, and as in Section \ref{sec-genincexc}, $c_4(G)$ is the number of cycles on $4$ vertices in $G$.
\begin{lemma} \label{lem-KwinsforP4-C4}
Fix $n$ and $d$ with $(d+1)|n$. For all $n$-vertex, $d$-regular graphs $G$ we have
$$
p_4(G) - c_4(G) \geq p_4\left(\frac{n}{d+1}K_{d+1}\right) - c_4\left(\frac{n}{d+1}K_{d+1}\right),
$$
with equality if and only if $G=\frac{n}{d+1}K_{d+1}$.
\end{lemma}

\medskip

\noindent {\em Proof}: Put an arbitrary ordering $<$ on the vertices of $G$. For each edge $e=uv$ of $G$ with $u<v$, denote by $A(e)$ those neighbors of $u$ that are not neighbors of $v$, by $B(e)$ those neighbors of $u$ that are also neighbors of $v$, and by $C(e)$ those neighbors of $v$ that are not neighbors of $u$. Set $k(e)=|B(e)|$ (so $0 \leq k(e) \leq d-1$ and $|A(e)|=|C(e)|=d-1-k(e)$). For each $x \in A(e)$ and $y \in C(e)$ there is a unique path on four vertices in $G$ that has $x$ and $y$ as end vertices and $e$ as middle edge, and the same is true for each $x \in A(e)$ and $y \in B(e)$, and for each $x \in C(e)$ and $y \in B(e)$. For each $x\neq y \in B(e)$ there are two such paths. In this way all paths on four vertices in $G$ with middle edge $e$ are counted. It follows that
\begin{eqnarray}
p_4(G) & = & \sum_{e \in E(G)} \left(|A(e)||C(e)|+|A(e)||B(e)|+|C(e)||B(e)|+2{|B(e)| \choose 2}\right) \nonumber \\
& = & \sum_{e \in E(G)} \left((d-1)^2-k(e)\right). \label{p4count}
\end{eqnarray}
Let $\ell(e)$ be the number of edges that go from $A(e)$ to $B(e)$, or from $B(e)$ to $C(e)$, or from $C(e)$ to $A(e)$, and let $m(e)$ be the number of edges that are inside $B(e)$. For each edge $e'$ counted by $\ell(e)$ there is a unique cycle on four vertices in $G$ passing through $e$ and $e'$. For each edge counted by $m(e)$ there are two such cycles. As we run over all edges of $G$, each cycle appears four times in this count, and so
$$
c_4(G) = \frac{1}{4} \sum_{e \in E(G)} \left(\ell(e)+2m(e)\right).
$$
Combining with (\ref{p4count}) we see that
$$
p_4(G)-c_4(G) = \frac{1}{4}\sum_{e \in E(G)} \left(4(d-1)^2 - \left(4k(e) + \ell(e) + 2m(e)\right) \right).
$$
For a fixed value of $k(e)$, the quantity $4k(e) + \ell(e) + 2m(e)$ is maximized when all edges are present between $A(e)$ and $C(e)$, $A(e)$ and $B(e)$, and $C(e)$ and $B(e)$, and so $\ell(e) = (d-1-k(e))^2+2k(e)(d-1-k(e))$, and all edges inside $B(e)$ are present, and so $m(e)={k(e) \choose 2}$, and the maximum value is
$(d-1)^2 + 3k(e)$. This in turn is maximized when $k(e)=d-1$. This maximum can be achieved for each $e$ only if for each edge in $G$, the end vertices of the edge share the same neighborhood, with that neighborhood inducing a complete graph. This occurs uniquely when $G=\frac{n}{d+1}K_{d+1}$.
\qed

\medskip

\noindent {\em Proof of Theorem \ref{thm-qstateWR}}: We begin by examining the term ${\rm hom}(G(S),H^c)$ in Lemma \ref{inq-incexc}.
Since $H^c$ is the disjoint union of complete bipartite subgraphs spanning no more than $\ell$ edges, we have $H^c=K_{r_1,s_1} \cup \ldots \cup K_{r_m,s_m}$ for some strictly positive integers $m$ and $r_i, s_i$ ($1 \leq i \leq m$) with $\sum_{i=1}^m r_is_i \leq \ell$. If $C$ is a bipartite component of $G(S)$ with partition classes having sizes $a$ and $b$ then
$$
{\rm hom}(C,H^c) = \sum_{i=1}^m \left(r_i^as_i^b+r_i^bs_i^a\right).
$$
If $C$ is not bipartite then ${\rm hom}(C,H^c)=0$. It follows that the contribution from $S$ to the sum in Lemma \ref{inq-incexc} will be non-zero only if $G(S)$ is bipartite.

From Lemma \ref{inq-incexc} and the calculation above we see that ${\rm hom}(G,H)$ is a polynomial in $q$ of degree $n$. The contribution to the polynomial from all $S$ that span at most $4$ vertices and do not span a path or a cycle on four vertices is independent of $G$. The contribution from each $S$ that spans a cycle on four vertices is
$$
\sum_{i=1}^m \left(2r_i^2s_i^2\right)q^{n-4}
$$
and the contribution from each $S$ that spans a path on four vertices is the negative of this. Using Lemma \ref{lem-KwinsforP4-C4} and reasoning in the same way as in the proof of Theorem \ref{thm_Hqell} we get that if $G \neq \frac{n}{d+1}K_{d+1}$ then
$$
{\rm hom}\left(\frac{n}{d+1}K_{d+1},H_q\right)-{\rm hom}(G,H_q) \geq Cq^{n-4} -Dq^{n-5}
$$
where $C=\sum_{i=1}^m(2r_i^2s_i^2)$ and $D$ is bounded above by the maximum value of
$$
\sum_{S \subseteq E(G)} {\rm hom}(G(S),H^c)
$$
over all $n$-vertex, $d$-regular graphs $G$, where the sum runs over all $S$ with $v(S)\geq 5$ that span a bipartite subgraph. Letting $C'$ be the minimum value of $C$ and $D'$ the maximum value of $D$ over all choices of $m$ and $r_i, s_i$ ($1 \leq i \leq m$), we see that as long as $q>D'/C'$ we have ${\rm hom}(\frac{n}{d+1}K_{d+1},H) > {\rm hom}(G,H)$, proving the theorem.

In the specific case of $H=H_q$ we have $m=r_1=s_1=1$ and so $C=2$, and since a bipartite subgraph of $G$ can have at most $n/2$ components (we are dealing with subgraphs without isolated vertices), we also have $D \leq \exp_2\{nd/2+n/2\}$. It follows that ${\rm hom}(\frac{n}{d+1}K_{d+1},H_q) > {\rm hom}(G,H_q)$ as long as $q > \exp_2\{nd/2+n/2-1\}$, as claimed just before the statement of Theorem \ref{thm-qstateWR}. \qed

\section{Comparing ${\rm hom}(K_{d,d},H)$ and ${\rm hom}(K_{d+1},H)$} \label{sec-recipe}

In this section we determine which of ${\rm hom}(K_{d,d},H)^{1/2d}$, ${\rm hom}(K_{d+1},H)^{1/(d+1)}$ is larger for all sufficiently large $d$, for each $H$. The determination will be in terms of some parameters that are fairly simply calculated from $H$ (given the ability to examine all subgraphs of $H$). Say that $H$ (without isolated vertices) is {\em of complete bipartite type} (respectively, {\em of complete type}, {\em of neutral type}) if there is some $d(H)>0$ such that for all $d \geq d(H)$ we have ${\rm hom}(K_{d+1},H)^{1/(d+1)} < {\rm hom}(K_{d,d},H)^{1/2d}$ (respectively, ${\rm hom}(K_{d,d},H)^{1/2d} < {\rm hom}(K_{d+1},H)^{1/(d+1)}$, ${\rm hom}(K_{d+1},H)^{1/(d+1)} = {\rm hom}(K_{d,d},H)^{1/2d}$).

We begin by considering loopless $H$. For such $H$ we have ${\rm hom}(K_{d+1},H)=0$ unless $H$ contains $K_{d+1}$ as a subgraph, which will not be the case for all sufficiently large $d$. On the other hand, since $H$ has at least one edge we will always have ${\rm hom}(K_{d,d},H)\geq 1$. So all loopless $H$ are of complete bipartite type, and from now on we assume that $H$ has at least one loop.

Our goal now is to obtain sharp estimates on ${\rm hom}(K_{d,d},H)$ and ${\rm hom}(K_{d+1},H)$. We begin with ${\rm hom}(K_{d,d},H)$.
Say that a pair $(A,B) \subseteq V(H)^2$ with $A, B \neq \emptyset$ is a {\em complete bipartite image} if for all $a \in A$ and $b \in B$ we have $ab \in E(H)$; the {\em size} of the pair is $|A||B|$. Define
\begin{eqnarray*}
\eta(H) & = & \mbox{the maximal size of a complete bipartite image in $H$, and} \\
m(H) & = & \mbox{the number of complete bipartite images of maximal size in $H$}.
\end{eqnarray*}
We have
$$
{\rm hom}(K_{d,d},H) = \sum \left\{S(d,A)S(d,B) : (A,B)~\mbox{a complete bipartite image} \right\}
$$
where $S(d,A)$ is the number of ways of coloring the set $\{1, \ldots, d\}$ with colors from $A$ in such a way that all of the colors of $A$ are used at least once. By inclusion-exclusion,
$$
S(d,A) = \sum_{i=0}^{|A|} (-1)^i{|A| \choose i}\left(|A|-i\right)^d
$$
and so
\begin{equation}  \label{inq-Kdd'}
{\rm hom}(K_{d,d},H) = m(H)\eta(H)^d + \sum_{i=0}^{\eta(H)-1} c_ii^d
\end{equation}
where the $c_i$'s are (not necessarily non-zero) constants.

Now we look at ${\rm hom}(K_{d+1},H)$. Say that a pair $(A,B) \subseteq V(H)^2$ is a {\em complete image} if $A$ and $B$ are disjoint, $A$ induces a complete looped graph, $B$ induces a complete unlooped graph, $A$ has at least one vertex, and $(A,B)$ is a complete bipartite image; the {\em primary size} of the pair is $|A|$ and the {\em secondary size} is $|B|$. Define
\begin{eqnarray*}
a(H) & = & \mbox{the maximal primary size of a complete image in $H$,} \\
b(H) & = & \mbox{the maximal secondary size of a complete} \\
&  & ~~~\mbox{image of maximal primary size, and} \\
n(H) & = & \mbox{the number of complete images in $H$ with} \\
&  & ~~~\mbox{primary size $a(H)$, secondary size $b(H)$}.
\end{eqnarray*}
Note that if $(A,B)$ is a complete image then $(A, A \cup B)$ is a complete bipartite image and so we always have
\begin{equation} \label{obsv}
a^2(H) + a(H)b(H) \leq \eta(H).
\end{equation}
This in particular means that $a^2(H) \leq \eta(H)$ always and that if $a^2(H) = \eta(H)$ then $b(H)=0$.

To obtain an analog of (\ref{inq-Kdd'}) for ${\rm hom}(K_{d+1},H)$ we begin by observing that for all large $d$
$$
{\rm hom}(K_{d+1},H) = \sum \left\{(d+1)_{|B|}S(d+1-|B|,A) : (A,B)~\mbox{a complete image} \right\},
$$
where $(x)_m=x(x-1)\ldots(x-(m-1))$.
Indeed, if $C$ is a subset of $V(H)$ that can be realized as $\{f(v): v \in V(K_{d+1})\}$ for some $f \in {\rm Hom}(K_{d+1},H)$, then $C$ must induce a complete graph (which for all large $d$ must have at least one loop). Writing $A(C)$ for the set of vertices of $C$ that have a loop and $B(C)$ for those that do not, each element of $B(C)$ must occur as $f(v)$ for exactly one $v \in V(K_{d+1})$. The factor $(d+1)_{|B|}$ above counts the number of ways of assigning unique preimages from $V(K_{d+1})$ to the vertices of $B$, and the factor $S(d+1-|B|,A)$ counts the number of ways of assigning preimages to the vertices of $A$. We therefore have
\begin{equation} \label{inq-Kdp1'}
{\rm hom}(K_{d+1},H)  = \frac{n(H)p_{a(H)}(H)}{a(H)^{b(H)-1}}a(H)^{d} + \sum_{i=0}^{a(H)-1} p_i i^{d}
\end{equation}
where the $p_i$'s are (not necessarily non-zero) polynomials in $d$ whose degrees and coefficients are constants, with in particular $p_{a(H)}(H)$ a monic polynomial of degree $b(H)$.

Comparing (\ref{inq-Kdd'}) and (\ref{inq-Kdp1'}), we see immediately that if $a^2(H)< \eta(H)$ then $H$ is of complete bipartite type. If not then by (\ref{obsv}) we have $a^2(H)=\eta(H)$ and $b(H)=0$ and so (\ref{inq-Kdd'}), (\ref{inq-Kdp1'}) yield
$$
{\rm hom}(K_{d,d},H)^{d+1} = m(H)^{d+1} (1+o(c_1^d))a(H)^{2d^2+2d}
$$
and
$$
{\rm hom}(K_{d+1},H)^{2d} = n(H)^{2d} (1+o(c_2^d))a(H)^{2d^2+2d}
$$
where the $c_i$'s (here and later) are positive constants (depending on $H$) strictly smaller than $1$. From this it immediately follows that if $m(H) >1$ and $m(H) \geq n(H)^2$ then $H$ is of complete bipartite type, while if $m(H) < n(H)^2$ then $H$ is of complete type.

It remains to study the case $a^2(H)=\eta(H)$ and $n(H)=m(H)=1$. In this case there is $A_0 \subseteq V(H)$ with $|A_0|=a(H)$ such that $(A_0,\emptyset)$ is the unique complete image of size $a(H)$ and $(A_0,A_0)$ is the unique complete bipartite image of size $\eta(H)$. Define
\begin{eqnarray*}
a'(H) & = &  \mbox{the largest primary size of a complete image $(A',B')$ with $A' \cup B' \not \subset A_0$,} \\
b'(H) & = &  \mbox{the largest secondary size of a such a complete image} \\
& & ~~~\mbox{of largest primary size, and} \\
n'(H) & = & \mbox{the number of such complete images with primary size $a'(H)$,} \\
& & ~~~\mbox{secondary size $b'(H)$.}
\end{eqnarray*}
If there are no such complete images then set $a'(H)=0$. Note that this is equivalent to $H$ being the disjoint union of a complete looped graph and a (perhaps empty) loopless graph.
Similarly define
\begin{eqnarray*}
\eta'(H) & = &  \mbox{the maximum size of a complete bipartite image $(A',B')$ with} \\
& & ~~~\mbox{at least one of $A', B' \not \subseteq A_0$, and} \\
m'(H) & = &  \mbox{the number of such complete bipartite images having size $\eta'(H)$.}
\end{eqnarray*}
If there are no such complete bipartite images then set $\eta'(H)=0$. Note that this is equivalent to $H$ being a complete looped graph.

If $\eta'(H)=0$ (which implies $a'(H)=0$) then trivially $H$ is of neutral type. If $\eta'(H) > 0$ and $a'(H)=0$ then ${\rm hom}(K_{d+1},H)^{1/(d+1)}=a(H)$ and ${\rm hom}(K_{d,d},H)^{1/2d} > \eta(H)^{1/2}=a(H)$
(because there is some contribution to ${\rm hom}(K_{d,d},H)$ from that part of $H$ that witnesses $\eta'(H)>0$ that is not counted among those colorings that take values only inside $A_0$), so $H$ is of complete bipartite type. So from now on we may assume that $\eta'(H)$, $m'(H)$, $a'(H)$ and $n'(H)$ are all strictly positive.

The dominant term of ${\rm hom}(K_{d+1},H)$ is now $a(H)^{d+1}$, which counts the number of homomorphisms which have as their image any subset of $A_0$, the unique complete looped subgraph of $H$ of size $a(H)$, and the sum of the remaining terms is
$$
\frac{n'(H)p'(H)}{a'(H)^{b'(H)}} (1+o(c_3^d))a'(H)^{d+1}
$$
where $p'(H)$ is a monic polynomial in $d$ of degree $b'(H)$ whose coefficients are constants. From this it follows that
\begin{equation} \label{kdp1-secorder}
\frac{{\rm hom}(K_{d+1},H)^{2d}}{a(H)^{2d^2+2d}} = 1 + (1+o(c_4^d))\frac{2dn'(H)p'(H)}{a(H)a'(H)^{b'(H)-1}}\left(\frac{a'(H)}{a(H)}\right)^d.
\end{equation}
Similarly the dominant term of ${\rm hom}(K_{d,d},H)$ is $\eta(H)^d=a(H)^{2d}$, which counts the number of homomorphisms which have as their image any subset of $A_0$, and the sum of the remaining terms is $m'(H)(1+o(c_5^d))\eta'(H)^d$. From this it follows that
\begin{equation} \label{kdd-secorder}
\frac{{\rm hom}(K_{d,d},H)^{d+1}}{a(H)^{2d^2+2d}} = 1 + (1+o(c_6^d))(d+1)m'(H)\left(\frac{\eta'(H)}{a^2(H)}\right)^d.
\end{equation}
Comparing (\ref{kdp1-secorder}) and (\ref{kdd-secorder}) we immediately get that if $a(H)a'(H) < \eta'(H)$ then $H$ is of complete bipartite type (note that this inequality can be interpreted as covering the case $\eta'(H)>0=a'(H)$), while if either $a(H)a'(H) > \eta'(H)$ or $a(H)a'(H) = \eta'(H)$ and $b'(H)>0$ then $H$ is of complete type.
If $a(H)a'(H) = \eta'(H)$ and $b'(H)=0$ then the right-hand side of (\ref{kdp1-secorder}) reduces to
$$
1 + (1+o(c_4^d))\frac{2dn'(H)a'(H)}{a(H)}C^d
$$
while the right-hand side of (\ref{kdd-secorder}) becomes $1+ (1+o(c_6^d))(d+1)m'(H)C^d$ (for some positive constant $C$, the same in both cases). We see from this that if $2n'(H)a'(H) > a(H)m'(H)$ then $H$ is of complete type, while if $2n'(H)a'(H) \leq a(H)m'(H)$ then $H$ is of complete bipartite type.

We summarize all of this in the following theorem. Here (as always) $H$ ranges over all simple, finite graphs, perhaps with loops, but without isolated vertices.
\begin{thm} \label{thm-char}
If $H$ is not a complete looped graph, then it is either of complete bipartite type or of complete type. The following are the conditions under which $H$ is of complete bipartite type:
\begin{enumerate}
\item $H$ is loopless.
\item $a^2(H) < \eta(H)$.
\item $a^2(H) = \eta(H)$, $m(H) >1$ and $m(H) \geq n(H)^2$.
\item $a^2(H) = \eta(H)$, $m(H) = n(H) = 1$ and $a(H)a'(H) < \eta'(H)$.
\item $a^2(H) = \eta(H)$, $m(H) = n(H) = 1$, $a(H)a'(H) = \eta'(H)$, \\$b'(H)=0$, $2n'(H)a'(H) \leq a(H)m'(H)$ and $m'(H)>0$.
\end{enumerate}
The following are the conditions under which $H$ is of complete type:
\begin{enumerate}
\item $a^2(H) = \eta(H)$ and $m(H) < n(H)^2$.
\item $a^2(H) = \eta(H)$, $m(H) = n(H) = 1$ and $a(H)a'(H) > \eta'(H)$.
\item $a^2(H) = \eta(H)$, $m(H) = n(H) = 1$, $a(H)a'(H) = \eta'(H)$ and $b'(H)>0$.
\item $a^2(H) = \eta(H)$, $m(H) = n(H) = 1$, $a(H)a'(H) = \eta'(H)$, \\$b'(H)=0$ and $2n'(H)a'(H) > a(H)m'(H)$.
\end{enumerate}
\end{thm}

We illustrate Theorem \ref{thm-char} with some specific examples, including all of the graphs considered so far in this paper.

The complete graph $K_q$ (which encodes proper $q$-colorings) has no loops and so is of complete bipartite type. More generally, the graph $H_q^\ell$ with $\ell>0$ (the subject of Theorem \ref{thm_Hqell}) is of complete bipartite type since it has $\eta(H)=q(q-\ell)$ and $a(H)=q-\ell$, so $a^2(H) < \eta(H)$.

The $k$-state hard-core constraint graph $H(k)$ ($k \geq 1$), is the graph on vertex set $\{0, \ldots, k\}$ with $ij \in E(H_k)$ if $i+j \leq k$. This graph occurs naturally in the study of multicast communication networks, and has been considered in \cite{MitraRamananSengupta} and \cite{GalvinMartinelliRamananTetali}. Note that for $k=1$ this is the same as $H_{\rm ind}$, the graph that encodes independent sets. When $k$ is odd, say $k=2\ell+1$, we have $\eta(H(k))=(\ell+1)(\ell+2)$ and $a(H(k))=\ell+1$ and so $H(2\ell+1)$ is of complete bipartite type. If $k$ is even, say $k=2\ell$, we have $\eta(H(k))=(\ell+1)^2$, $m(H(k))=1$, $a(H(k))=\ell+1$ and $n(H(k))=1$ and so we have to look at the primed parameters. We have $a'(H(k))=\ell$ and $\eta'(H(k))=(\ell+1)^2-1$ and so $a(H)a'(H)< \eta'(H)$ and $H(2\ell)$ is also of complete bipartite type.

The completely disconnected fully looped graph $E_k^o$ ($k \geq 1$) is the graph consisting of $k$ loops and no other edges. For $k=1$, $E_k^o$ is trivially of neutral type. For $k>1$ we have $\eta(H)=1$, $a(H)=1$, $m(H)=k$ and $n(H)=k$, so that for all $k>1$ it is of complete type. This was the example that was first pointed out by Cutler and Radcliffe (personal communications), showing that Galvin and Tetali's original conjecture \cite{GalvinTetali-weighted} (concerning the validity of (\ref{GTbound}) for non-bipartite $G$) was false.

The complete looped path $P_k^o$ ($k \geq 1$) is the path on $k$ vertices with all vertices looped. For $k=1$ and $2$ this coincides with the complete looped graph on $k$ vertices, and $P_k^o$ in these cases is trivially of neutral type. For $k \geq 3$ we have $\eta(H)=4$, $a(H)=2$, and $m(H)=n(H)=k-1$, so $P_k^o$ in these cases is of complete type. Thus one of the simplest non-trivial examples of a graph of complete type is the fully looped path on three vertices, also known as the Widom-Rowlinson graph $H_{\rm WR}$. Similarly the graph $H_q$ ($q \geq 3$) (the subject of Theorem \ref{thm-qstateWR}) has $\eta(H)=(q-1)^2$, $a(H)=q-1$, and $m(H)=n(H)=2$, so $H_q$ is also of complete type.

\section{Proofs of Theorems \ref{thm-asymptotic} and \ref{thm-WR}} \label{sec-asymptotics}

We make use of the following result, a special case of a theorem of Madiman and Tetali \cite[Theorem III]{MadimanTetali}. The result was originally obtained for $H=H_{\rm ind}$ by Kahn (unpublished), and the present author first observed the generalization to all $H$.
\begin{thm} \label{thm-MT}
Let $G$ be a $d$-regular graph and $H$ an arbitrary graph. Let $<$ be an arbitrary ordering of the vertices of $G$, and let $p(v)$ denote the number of neighbors $u$ of $v$ with $u < v$. We have
$$
{\rm hom}(G,H) \leq \prod_{v \in V(G)~\!:~\!p(v)>0} {\rm hom}(K_{p(v),p(v)}, H)^{\frac{1}{d}}.
$$
\end{thm}
(The exclusion of those $v$ with $p(v)=0$ does not appear in \cite{MadimanTetali}, but the result stated above follows immediately since ${\rm hom}(K_{0,0},H)=1$ for all $H$.)

By (\ref{inq-Kdd'}) we have that for $i>0$,
$$
{\rm hom}(K_{i,i}, H) \leq m(H)\eta(H)^{i}\left(1+ C(H)(1-1/\eta(H))^{i}\right)
$$
for some constant $C(H)>0$ and so, using $\sum_{v \in V~\!:~\!p(v)>0} p(v) = |E(G)| = nd/2$ and $\log(1+x)\leq x$, we get from Theorem \ref{thm-MT} that
\begin{equation} \label{inq-MT2}
{\rm hom}(G,H) \leq \eta(H)^{\frac{n}{2}} m(H)^{\frac{n-p_0n}{d}} \exp\left\{\frac{C(H)n}{d}\sum_{i=1}^d p_if(i)\right\},
\end{equation}
where $p_i$ is the proportion of $v \in V(G)$ with $p(v)=i$, $n=|V(G)|$, and $f(i)=(1-1/\eta(H))^i$.

Our approach now is to choose an ordering of the vertices of $G$ that makes the right-hand side of (\ref{inq-MT2}) small. We begin by fixing an independent set $I$ of maximal size, say $\alpha n$, and putting the vertices of $I$ at the beginning of the ordering (in some arbitrary order). We then extend to a total ordering by choosing a (uniform) random ordering of the vertices outside $I$. Since $I$ is maximal every $v \not \in I$ has a neighbor in $I$ and so $p_0=\alpha$.

Let $A_i$ ($i=1, \ldots, d$) be the set of vertices outside of $I$ that have exactly $i$ neighbors in $I$, and let $a_in=|A_i|$. For each $v \in A_i$ we have (by a simple symmetry argument) that $\Pr(p(v)=j)=1/(d-(i-1))$ for $j=i, \ldots, d$ (and $\Pr(p(v)=j)=0$ otherwise). It follows that the expected value of $p_i$ is $\sum_{j=1}^i a_j/(d-(j-1))$, and so by linearity of expectation the expected value of $\sum_{i=1}^d p_if(i)$ is
$$
\sum_{i=1}^d \sum_{j=1}^i \frac{a_jf(i)}{d-(j-1)} = \sum_{j=1}^d \frac{a_j}{d-(j-1)}\sum_{i=j}^d f(i).
$$
Choose an ordering for which $\sum_{i=1}^d p_if(i)$ is at most this much. Since $f(i) = (1-1/\eta(H))^i$ we have $\sum_{i=j}^d f(i) \leq O((1-1/\eta(H))^j)$. For sufficiently large $d$, the maximum of $(1-1/\eta(H))^j/(d-(j-1))$ as $j$ ranges from $1$ to $d$ is achieved at $j=1$. Using this in (\ref{inq-MT2}) together with the fact that $\sum_{i=1}^d a_i \leq 1$ leads to
\begin{equation} \label{inq-MT4}
{\rm hom}(G,H)^\frac{1}{n} \leq \eta(H)^\frac{1}{2} m(H)^\frac{1-\alpha}{d} \exp\left\{\frac{C'(H)}{d^2}\right\} = (1+o(1))\eta(H)^{\frac{1}{2}}
\end{equation}
for some constant $C'(H)>0$.

We now compare the bound in (\ref{inq-MT4}) to asymptotic bounds on ${\rm hom}(K_{d,d},H)^{1/2d}$ (or ${\rm hom}(K_{d+1},H)^{1/(d+1)}$, as appropriate), derived from identities in Section \ref{sec-asymptotics}. We exclude two trivial cases that have been dealt with completely earlier, namely $H$ bipartite and $H$ completely looped.

If $H$ is of complete bipartite type with $m(H)>1$, then by (\ref{inq-Kdd'}) we have
\begin{equation} \label{asy1}
{\rm hom}(K_{d,d},H)^{1/2d} = \eta(H)^\frac{1}{2} m(H)^{\frac{1+o(1)}{2d}},
\end{equation}
so the upper bound in (\ref{inq-MT4}) is correct up to the exponent of $m(H)$ in the second term.

\medskip

At this point we draw attention to perhaps the simplest instance of a graph of complete bipartite type for which Conjecture \ref{conj-GTcorr} remains open, namely $H=K_3$ (the case of proper $3$-coloring). The conjecture is trivially true for those $G$ without an independent set of size at least $n/3$ (for which ${\rm hom}(G,K_3)=0$). For the remaining $G$, the upper bound in (\ref{inq-MT4}) becomes
$$
{\rm hom}(G,K_3)^\frac{1}{n} \leq 2^\frac{1}{2} 6^\frac{2+o(1)}{3d},
$$
whereas
$$
{\rm hom}(K_{d,d},K_3)^{1/2d} = 2^\frac{1}{2} 6^\frac{1+o(1)}{2d}.
$$
It would be very nice to be able to close the gap between the $6^{2/3}$ and $6^{1/2}$ here.

\medskip

If $H$ is of complete bipartite type with $m(H)=1$, then by (\ref{inq-Kdd'}) we have
\begin{equation} \label{asy2}
{\rm hom}(K_{d,d},H)^{1/2d} = \eta(H)^\frac{1}{2} \exp\left\{o(c_7^d)\right\},
\end{equation}
so the upper bound in (\ref{inq-MT4}) is correct only in the first term.

We now examine those $H$ of complete type. First we look at the case $a^2(H) = \eta(H)$ and $m(H) < n(H)^2$. Here $m(H)\geq n(H)$ (since every complete looped graph $A$ of size $a(H)$ gives rise to a pair $(A,A)$ that gets counted in $m(H)$), and so $m(H)>1$. Also in this case $b(H)=0$ and $a(H)=\eta(H)^{1/2}$ and so from (\ref{inq-Kdp1'}) we get
\begin{equation} \label{asy3}
{\rm hom}(K_{d+1},H)^{1/(d+1)} = \eta(H)^\frac{1}{2}n(H)^{\frac{1+o(1)}{d}}.
\end{equation}
If $m(H)>n(H)$ then the upper bound in (\ref{inq-MT4}) is correct up to the exponent of $m(H)$ in the second term. If, however $m(H)=n(H)$, then the upper bound in (\ref{inq-MT4}) has the correct exponent of $m(H)$, and it is in later terms that we first see disagreement. In this case, a more precise comparison between (\ref{inq-MT4}) and ${\rm hom}(K_{d+1},H)$ is in order. From (\ref{inq-Kdp1'}) (and using $b(H)=0$) we have
\begin{equation} \label{inq-Kdp1''}
{\rm hom}(K_{d+1},H)^{1/(d+1)}  = \eta(H)^{\frac{1}{2}}n(H)^\frac{1}{d+1} \exp\left\{o(c_8^d)\right\}
\end{equation}
for large $d$. Defining $\gamma$ by $\alpha=1/(d+1)+\gamma$, (\ref{inq-MT4}) becomes
\begin{equation} \label{inq-MT3}
{\rm hom}(G,H)^\frac{1}{n} \leq \eta(H)^\frac{1}{2} m(H)^\frac{1}{d+1} \exp\left\{-\frac{C''(H)\gamma}{d}+\frac{C'(H)}{d^2}\right\}
\end{equation}
where $C''(H)>0$ is a constant. Comparing (\ref{inq-MT3}) and (\ref{inq-Kdp1''}), we see that for $\gamma > \Omega(1/d)$ (with the constant depending on $H$) we have the exact bound
$$
{\rm hom}(G,H) \leq {\rm hom}(K_{d+1},H)^{\frac{n}{d+1}}.
$$
We have proven the following generalization of Theorem \ref{thm-WR}.
\begin{thm} \label{thm-WR-gen}
Let $H$ be of complete type with $m(H)=n(H)>1$. There is a constant $C>0$ such that for all sufficiently large $d$, we have the following. If $G$ is an $n$-vertex, $d$-regular graph with an independent set of size at least $Cn/d$, then
$$
{\rm hom}(G,H) < {\rm hom}(K_{d+1},H)^\frac{n}{d+1}.
$$
\end{thm}
Note that $H_{\rm WR}$, and more generally $H_q$ of Theorem \ref{thm-qstateWR} for all $q\geq 3$, satisfies the conditions of Theorem \ref{thm-WR-gen}, as does for example the complete looped path on $k \geq 3$ vertices.

There remains the case of $H$ of complete type with $m(H) = n(H) = 1$. In this case, we have
\begin{equation} \label{asy4}
{\rm hom}(K_{d+1},H)^{1/(d+1)} = \eta(H)^\frac{1}{2} \exp\left\{o(c_9^d)\right\},
\end{equation}
so the upper bound in (\ref{inq-MT4}) is correct only in the first term.

From (\ref{asy1}), (\ref{asy2}), (\ref{asy3}) and (\ref{asy4}) we see that in all cases
$$
\max\left\{{\rm hom}(K_{d,d},H)^{\frac{1}{2d}}, {\rm hom}(K_{d+1},H)^{\frac{1}{d+1}}\right\} = (1+o(1))\eta(H)^\frac{1}{2},
$$
which together with (\ref{inq-MT4}) gives Theorem \ref{thm-asymptotic}.

\bigskip

\noindent {\bf Acknowledgement}: We are grateful to Felix Lazebnik for suggesting the use of Whitney's broken circuit theorem.

\end{document}